\documentclass{gtart_h}

\def\ifplaintex{\expandafter\ifx\csname documentclass\endcsname\relax}

\def\gtp{{\mathsurround=0pt\it $\cal G\mskip-2mu$eometry \&\ 
$\cal T\!\!$opology $\cal P\!$ublications}}  

\def\recd{{\small Received:\qua\receiveddate\ifx\reviseddate\relax
\else\qquad Revised:\qua\reviseddate\fi\par}} 

\def\lognumber#1{\def\thelognumber{#1}}
\def\volumenumber#1{\def\thevolumenumber{#1}}
\def\volumeyear#1{\def\thevolumeyear{#1}}
\def\papernumber#1{\def\thepapernumber{#1}}
\def\pagenumbers#1#2{\def\startpage{#1}\def\finishpage{#2}}
\def\published#1{\def\publishdate{#1}}

\def\received#1{\def\receiveddate{#1}}
\def\revised#1{\def\reviseddate{#1}}
\def\accepted#1{\def\accepteddate{#1}}

\def\asciiaddress#1{\def\theasciiaddress{#1}}
\def\asciiemail#1{\def\theasciiemail{#1}}

\long\def\asciiabstract#1{\long\def\theasciiabstract{#1}}

\let\\\par\let\thelognumber\relax\let\thevolumenumber\relax
\let\thepapernumber\relax\let\thevolumeyear\relax\let\startpage\relax
\let\finishpage\relax\let\publishdate\relax\let\receiveddate\relax
\let\reviseddate\relax\let\accepteddate\relax\let\theasciititle\relax
\let\theasciiauthors\relax\let\theasciiaddress\relax
\let\theasciiabstract\relax

\let\theasciiemail\relax

\ifplaintex
\font\logobig=cmssbx10 scaled 3836
\font\logomed=cmssbx10 scaled 2557
\else
\font\logobig=cmssbx10 scaled 4200
\font\logomed=cmssbx10 scaled 2800
\fi

\long\def\makeagttitle{   
\count0=\startpage
\agt\hfill      
\hbox to 45truept{\vbox to 0pt{\vglue -13truept{\logomed A\kern -.37em{\logobig 
T}\kern -.38em G}\vss}\hss}
\break
{\small Volume \thevolumenumber\ (\thevolumeyear)
\startpage--\finishpage\nl
Published: \publishdate}

\vglue .25truein

{\parskip=0pt\leftskip 0pt plus
1fil\def\\{\par\smallskip}{\Large\bf\thetitle}\par\medskip} \vglue
0.05truein

{\parskip=0pt\leftskip 0pt plus 1fil\def\\{\par}{\sc\theauthors}
\par\medskip}%
 
\vglue 0.03truein 

{\small\leftskip 25truept\rightskip 25truept{\bf Abstract}\stdspace\theabstract

{\bf AMS Classification}\stdspace\theprimaryclass
\ifx\thesecondaryclass\relax\else; \thesecondaryclass\fi\par
{\bf Keywords}\stdspace \thekeywords\par}\vglue 7truept

}   

\ifplaintex
\font\phead=cmsl9 scaled 950
\font\pnum=cmbx10 scaled 913
\font\pfoot=cmsl9 scaled 950
\headline{\vbox to 0pt{\vskip -4.5mm\line{\small\phead\ifnum
\count0=\startpage ISSN 1472-2739 (on-line) 1472-2747 (printed)
\hfill {\pnum\folio}\else\ifodd\count0\def\\{ }%
\ifx\theshorttitle\relax\thetitle\else\theshorttitle\fi\hfill{\pnum\folio}
\else\def\\{ and }{\pnum\folio}\hfill\ifx\theshortauthors\relax\theauthors
\else\theshortauthors\fi\fi\fi}\vss}}
\footline{\vbox to 0pt{\vglue 0mm\line{\small\pfoot\ifnum\count0=\startpage
\copyright\ \gtp\hfill\else
\agt, Volume \thevolumenumber\ (\thevolumeyear)\hfill\fi}\vss}}
\else
\font=cmsl9 scaled 1050
\font\lnum=cmbx10 
\font=cmsl9 scaled 1050
\makeatletter
\def\@oddhead{{\small\lhead\ifnum\count0=\startpage ISSN 1472-2739 
(on-line) 1472-2747 (printed)\hfill {\lnum\number\count0}\else\ifodd\count0
\def\\{ }\ifx\theshorttitle\relax \thetitle \else\theshorttitle\fi\hfill
{\lnum\number\count0}\else\def\\{ and }{\lnum\number\count0}
\hfill\ifx\theshortauthors\relax 
\theauthors\else\theshortauthors\fi\fi\fi}}\def\@evenhead{\@oddhead}
\def\@oddfoot{\small\lfoot\ifnum\count0=\startpage\copyright\ \gtp\hfill\else
\agt, Volume \thevolumenumber\ (\thevolumeyear)\hfill\fi}
\def\@evenfoot{\@oddfoot}
\makeatother
\fi
\let\maketitlepage\makeagttitle
\let\makeshorttitle\maketitlepage
\let\maketitle\maketitlepage

\newwrite\gtoutfile
\long\gdef\makeheadfile{  
{\def\\{, }\def\s{ }
\immediate\openout\gtoutfile head.xxx
\immediate\write\gtoutfile{Proxy-for: \ifx\theasciiauthors\relax
\theauthors\else\theasciiauthors\fi\s<\ifx\theasciiemail\relax\theemail\else\theasciiemail\fi>}
\immediate\write\gtoutfile{\noexpand\\}
\immediate\write\gtoutfile{Authors: \ifx\theasciiauthors\relax
\theauthors\else\theasciiauthors\fi}
{\def\\{ }\immediate\write\gtoutfile{Title: \ifx\theasciititle\relax
\thetitle\else\theasciititle\fi}}
\immediate\write\gtoutfile{Subj-class: GT or SG, GR etc}
\immediate\write\gtoutfile{MSC-class: \theprimaryclass\ifx\thesecondaryclass\relax\else, \thesecondaryclass\fi}
\immediate\write\gtoutfile{Journal-ref: Algebr. Geom. Topol. \thevolumenumber\s
(\thevolumeyear) \startpage-\finishpage}
\immediate\write\gtoutfile{Comments: Published by Algebraic and
Geometric Topology at}
\immediate\write\gtoutfile{\s\s\s  http://www.maths.warwick.ac.uk/agt/AGTVol\thevolumenumber/agt-\thevolumenumber-\thepapernumber.abs.html}
\immediate\write\gtoutfile{\noexpand\\}
\immediate\write\gtoutfile{}
\ifx\theasciiabstract\relax
\immediate\write\gtoutfile{\theabstract}\else
\immediate\write\gtoutfile{\theasciiabstract}\fi
\immediate\write\gtoutfile{}
\immediate\write\gtoutfile{\noexpand\\}
\immediate\write\gtoutfile{}
\immediate\closeout\gtoutfile}}  

\def\maketitlepage{\makeagttitle\makeheadfile}
\let\makeshorttitle\maketitlepage
\let\maketitle\maketitlepage

\lognumber{57}
\volumenumber{5}
\volumeyear{2005}
\papernumber{57}
\pagenumbers{1451}{1469}
\received{25 March 2005} 
\revised{4 August 2005}
\accepted{14 September 2005}
\published{30 October 2005}

\usepackage{amssymb,amsmath,graphicx} 

\def\fref#1{\hyperlink{#1anchor}{\ref*{#1}}}
\def\figref#1{\hyperlink{#1anchor}{Figure~\ref*{#1}}}
\def\anchor#1{\noindent\hypertarget{#1anchor}{\smash{$\phantom{99}$}}}

\def\g{\gamma}
\def\a{\alpha}
\def\b{\beta}

\newtheorem{thm}{Theorem}[section]  
\newtheorem{lem}[thm]{Lemma}        
\newtheorem{cor}[thm]{Corollary}    
\newtheorem{prop}[thm]{Proposition}  

\theoremstyle{definition}
\newtheorem{defn}[thm]{Definition}

\begin{document}

\title{Hyperbolic covering knots}                    
\author{Daniel S. Silver\\Wilbur Whitten}
\address {Department of Mathematics, University of South 
Alabama\\Mobile AL 36688, USA}
\secondaddress{1620 Cottontown Road, Forest VA 24551, USA}                  
\asciiaddress{Department of Mathematics, University of South Alabama\\Mobile AL 36688, USA\\and\\1620 Cottontown Road, Forest VA 24551, USA}                  
\gtemail{\mailto{silver@jaguar1.usouthal.edu}, \mailto{bjwcw@aol.com}}  
\asciiemail{silver@jaguar1.usouthal.edu, bjwcw@aol.com }   

\begin{abstract}  
Given any knot $k$, there exists a hyperbolic knot $\tilde k$ with
arbitrarily large volume such that the knot group $\pi k$ is a
quotient of $\pi \tilde k$ by a map that sends meridian to meridian
and longitude to longitude. The knot $\tilde k$ can be chosen to be
ribbon concordant to $k$ and also to have the same Alexander invariant
as $k$.
\end{abstract}

\asciiabstract{%
Given any knot k, there exists a hyperbolic knot tilde k with
arbitrarily large volume such that the knot group pi k is a
quotient of pi tilde k by a map that sends meridian to meridian
and longitude to longitude. The knot tilde k can be chosen to be
ribbon concordant to k and also to have the same Alexander invariant
as k.}

\maketitle

\section{Introduction}
\label{intro}
 The classical problem of topology to find all homotopy classes of maps $M \to N$  between given complexes $M$ and $N$ has been variously expanded in recent years for the case in which $M$ and $N$ are manifolds of the same dimension; for an overview, see \cite{wang}. In the spirit of this expanded viewpoint as applied to knot theory, the authors in \cite{SW}  showed that given any knot $k$, there exists infinitely many prime knots $\tilde k$ admitting an epimorphism of knot groups $\pi {\tilde k} \to \pi k$ sending a meridian-longitude pair for $\tilde k$ to a meridian-longitude pair for $k$. We make use of this result, and go further, proving that the knots $\tilde k$ can in fact be chosen to be hyperbolic with arbitrarily large volumes (see Theorem 2.2).  The knots $\tilde k$ that we construct are ribbon concordant to $k$, and have the same Alexander invariant as $k$; in particular, they have the same Alexander polynomial. 

E. Kalfagianni showed in \cite{kalfagianni} that given any positive integer $n$, there exists a hyperbolic knot with trivial Alexander polynomial, trivial finite type invariants of orders $\le n$ and volume greater than $n$. Our result can be seen as a partial generalization.  

We are grateful to Abhijit Champanerkar, Tim Cochran and Danny Ruberman for helpful discussions. The first author is partially supported by NSF grant
DMS-0304971. 

{\bf Note added in proof}\qua Professor A. Kawauchi has informed the authors that many of the results of this paper can be found in 
\cite{Kaw1} or \cite{Kaw2}.


\section{Statement of results} 
\label{results}
We denote the group $\pi_1(S^3\setminus {\rm Int}(V), *)$ of a knot $k\subset S^3$ by $\pi k$. Here $V\cong k \times D^2$  is a tubular neighborhood of $k$, and $*$ is a basepoint chosen on the boundary $\partial V\cong k \times S^1$. An essential simple closed curve in $\partial V$ that is contractible in $V$ is called a {\it meridian}, and it is denoted by $m$. 
An essential simple closed curve $l \subset \partial V$ that is nullhomologous in $S^3\setminus {\rm Int}(V)$ is called a {\it longitude}. Once $k$ is oriented, both $m$ and $l$ acquire induced orientations. The inclusion map $ \partial V \hookrightarrow S^3 \setminus {\rm Int}(V)$ induces an  injection of fundamental groups. Its image
is the subgroup $\langle m, l\rangle$ generated by  $m$ and $l$. 

Let $k_i  (i= 1,2)$ be knots with meridian-longitude pairs $m_i,
l_i$.  

\begin{defn} A homomorphism $\phi:
\pi k_1 \to \pi k_2$ {\it preserves peripheral structure} if the
image of $\langle m_1, l_1\rangle$ is conjugate to a subgroup of
$\langle m_2,  l_2\rangle$. When $\phi$ is an epimorphism, we write $k_1 \succeq k_2$. 
\end{defn}

The relation $\succeq$ is a partial order \cite{SW}. After an appropriate choice of orientation, we can assume that $\phi(m_1)= m_2 l_2^p$ and $\phi(l_1)=m_2^ql_2^r$, for some integers $p, q,r$. Since $m_2^q l_2^r$ must be in $(\pi k_2)'' \cap Z(m_2)$ \cite{JL}, we have $q=0$. Furthermore, since the normal subgroup of $\pi {k_2}$ generated by $m_2l_2^p$ is all of $\pi k_2$, Corollary 2 of \cite{CGLS} implies that $p\in \{0,1,-1\}$; in fact the recent proof that every nontrivial knot satisfies Property P \cite{KM} implies that $p =0$.  Hence $\phi(m_1)= m_2$ and $\phi(l_1)=l_2^r$.  When $r=1$, we write $k_2 \succeq_1 k_1$. In \cite{SW} we showed that $k_1 \succeq_1 k_2$ implies $k_1 \succeq k_2$ but not conversely. 

A {\it ribbon concordance} from a knot $k_1$ to another knot $k_0$ is a smooth concordance $C\subset {\mathbb S}^3\times I$ with $C \cap {\mathbb S}^3 \times \{i\} = 
k_i (i = 0, 1)$, and such that the restriction to $C$ of the projection ${\mathbb S}^3 \times I\to I$ is a Morse function with no local maxima.  Visualizing such a concordance by cross-sections, we see a sequence of saddle points (called fusions) and local minima (the result of shrinking to points unknotted, unlinked components). We do not see any local maxima. 

The notion of ribbon concordance was introduced by C. Gordon \cite{gordon}, who wrote $k_1 \ge k_0$ if there is a ribbon concordance from $k_1$ to $k_0$. The term was motivated by the fact that a knot $k$ is ribbon concordant to the trivial knot if it bounds an immersed disk in ${\mathbb S}^3$ with only ribbon singularities. Gordon conjectured that  $\ge$ is a partial order. The conjecture remains open.  It is immediate from \cite{miyazaki} that ribbon concordance does not imply $\succeq$, nor does $\succeq$ imply ribbon concordance.

\begin{thm}\label{main}  Let $k$ be a knot. There exists a hypberbolic knot $\tilde k$ with the following properties. 
\begin{itemize}
\item[\rm(i)]  $\tilde k \succeq_1 k$; 
\item[\rm(ii)] The Alexander invariants of $\tilde k$ and $k$ are isomorphic; 
\item[\rm(iii)] $\tilde k$ has arbitrarily large volume;
\item[\rm(iv)] $\tilde k$ is ribbon concordant to $k$.
\end{itemize}
\end{thm} 

The $4$-{\it ball genus} of a knot $k\subset {\mathbb S}^3 = \partial B^4$ is the minimum genus of any properly embedded surface $F \subset B^4$ bounding $k$.

\begin{cor} Every Alexander polynomial is realized by hyperbolic knots with arbitrarily large volume and arbitrarily large 4-ball genus. 
\end{cor}

Corollary 2.3 is proven using results of J. Rasmussen \cite{rasmussen} and C. Livingston \cite{livingston}.  The statement of Corollary 2.3 was shown earlier by  A. Stoimenow using more combinatorial methods.


\section{Proof of Theorem 2.2} The idea for the proof Theorem 2.2 was suggested by \cite{miyazaki}. The rough idea is as follows. First, we invoke \cite{SW} so that we may assume without loss of generality that $k$ is prime. Having chosen a diagram for $k$ with a minimal number of crossings, we introduce a carefully devised unknot (called a ``staple") into a small neighborhood of each crossing. The greater part of the proof is devoted to showing that the resulting link is hyperbolic. Finally, we perform $1/q$ surgery on each of the staples. Thurston's hyperbolic surgery theorem implies that the resulting knots $\tilde k$ will be hyperbolic provided that the values of $q$ are sufficiently large. The special form of the staples ensures that $\tilde k$ has the same abelian invariants as $k$.

The main result of \cite{SW} implies that there exists a prime knot $\tilde k$ such that $\tilde k \succeq_1 k$. In fact there are infinitely many. Hence we can assume without any loss of generality that $k$ is prime.

\begin{figure}[ht!]\anchor{hck1}
\begin{center}
\includegraphics[width=1.6 in]{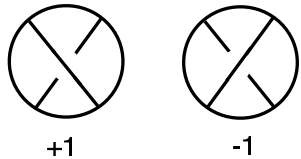}
\vspace{-4mm}\caption{Tangle $(B_i, t_i)$}
\label{hck1}
\end{center}
\end{figure}

Take a regular projection of $k$ with a minimal number $m$ of crossings. We may assume that $k$ lies in the projection plane except near the crossings. Number the crossings $i=1, \ldots, m$, and for each $i$, let $B_i$ be a $3$-ball that meets $k$ in two subarcs $t_{i_1}$ and $t_{i_2}$ that form the $i$th crossing. Thus each $(B_i, t_{i_1}\cup t_{i_2})$, abbreviated by $(B_i, t_i)$, is either the tangle $+1$ or $-1$, depending on the crossing (\figref{hck1}). We also assume that each $B_i$ meets the projection plane in an equatorial disk, and that $B_i \cap B_j =\emptyset$ when $i\ne j$. We assume that the balls $B_i$ 
are chosen so that $k \setminus t_1\cup \cdots \cup t_m$ is in the projection plane. 

\begin{figure}[ht!]\anchor{hck2}
\begin{center}
\includegraphics[width=3.2 in]{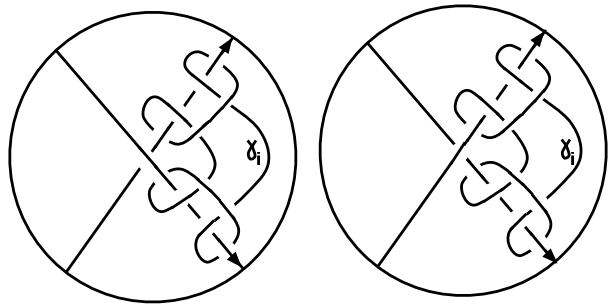}
\caption{Tangle $(B_i, t_i, \g_i)$}
\label{hck2}
\end{center}
\end{figure}

Next we insert an unknot $\g_i$ in the interior of each $B_i \setminus k_i$, as in \figref{hck2}. We refer to $\g_i$ as a {\it staple}. We orient $k$ in order to make the location of each staple specific. Note that $(B_i, t_i, \g_i)$ is homeomorphic to $(B_j, t_j, \g_j)$, for each $i$ and $j$. 

The proof of Theorem 2.2 proceeds by a sequence of lemmas.

\begin{lem} The link $L=k\cup \g_1\cup \cdots \cup \g_m$ is unsplittable. 
\end{lem}

\begin{proof} By construction, the sublink $\g_1\cup \cdots \cup \g_m$ is trivial. It suffices to show that $k\cup \g_i$ is unsplittable, for each $i$. 

It is convenient to have another view of $(B_i, t_i, \g_i)$, obtained
in the style of Montisenos by stretching $\partial B_i$ into an
``arc," as in \figref{hck3}a. \figref{hck3}b gives a view of the
2-fold cover of $B_i\setminus \g_i$ branched over $t_i$. It is a solid
torus $V_i$ minus the 2-component link $\tilde \g_i = \tilde
\g_{i_1}\cup \tilde \g_{i_2}$. The program Snap shows that $\tilde
\g_i$ is a hyperbolic link in $V_i$; that is, ${\rm Int}(V_i\setminus
\tilde \g_i)$ is a hyperbolic $3$-manifold.

\begin{figure}[ht!]\anchor{hck3}
\begin{center}
\includegraphics[width=3.8 in]{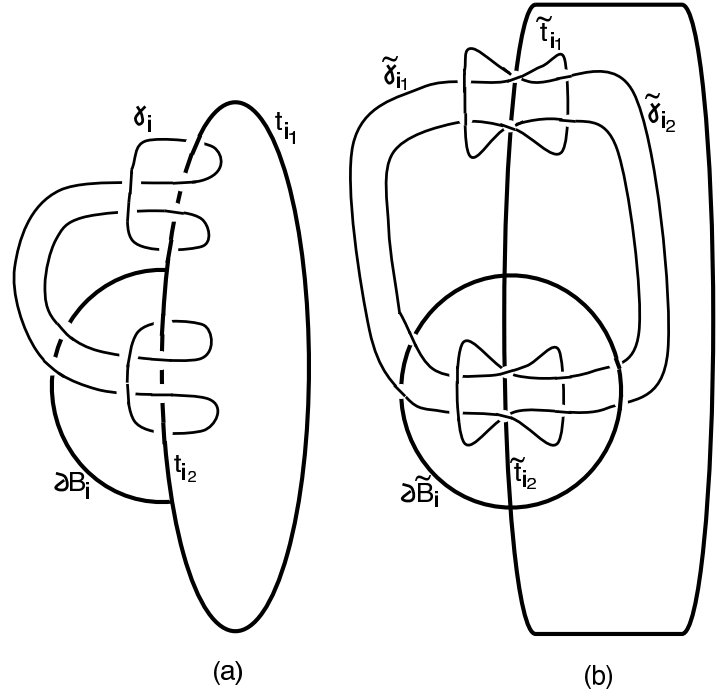}
\caption{(a)\qua Tangle $(B_i, t_i, \g_i)$ \quad (b)\qua 2-Fold branched cover}
\label{hck3}
\end{center}
\end{figure}

If $k\cup \g_i$ is splittable, then there exists a 2-sphere $S$ bounding a pair of $3$-balls, one containing $k$, the other, which we call $A$, containing $\g_i$. Since each of $B_i$ and $A$ contains $\g_i$, their interiors intersect. Clearly $B_i$ is not a subset of $A$, as $B_i$ contains two subarcs of $k$. Therefore if $A$ is not a subset of $B_i$, we can assume that $S\cap \partial B_i$ is a finite collection of pairwise disjoint simple closed curves. Let $\a$ be one of the curves that is innermost in $S$. 

If $\a$ bounds a disk $D$ in $S\cap{\rm cl}(S^3\setminus B_i)$, then it also bounds a disk $D'$ in $\partial B_i$ that is in $A$, and since $D'\cap k =\emptyset$, the sphere $D\cup D'$ bounds a $3$-ball not containing $k\cup \g_i$. Isotoping $D$ through the ball, we can remove $\a$ without moving $k\cup \g_i$. 

If, on the other hand, $\a$ bounds a disk $D\subset B_i$, then $\a$ also bounds a disk $D'\subset \partial B_i$ that contains no points of $t_i\cap \partial B_i$, since otherwise either $D\cap t_i\ne \emptyset$ or else $D$ lifts to a pair of meridianal disks of $V_i$ neither of which meets $\tilde \g_{i_1}\cup \tilde \g_{i_2}$. But $D\cap t_i = \emptyset$ by construction, and $\tilde \g_{i_1}\cup \tilde \g_{i_2}$ is essential in the $2$-fold cover of $B_i$ branched over $t_i$. Hence $D\cup D'$ bounds a $3$-ball $A' \subset B_i\setminus t_i$. If $\g_i\subset A'$, then we push $D'$ slightly into $B_i$ and replace $S$ by $D\cup D'$. If $\g_i$ is not a subset of $A'$, then we push $D$ through $A'$ into ${\rm cl}(S^3\setminus B_i)$, and thereby eliminate $\a$. 

Inductively, we remove all curves of $S\cap \partial B_i$, and assume henceforth that $S$ and hence $A$ are contained in the interior of $B_i$. However, the lift of $S$ to the $2$-fold cover $V_i\setminus \tilde \g_i$ of $B_i\setminus \g_i$ branched over $t_i$ is a pair of $2$-spheres, each of which splits $\tilde \g_i = \tilde \g_{i_1}\cup \tilde \g_{i_2}$. Since ${\rm Int}(V_i\setminus \tilde g_i)$ is hyperbolic and hence irreducible, this is impossible. Therefore, $k\cup \g_i$ is unsplittable. \end{proof}

\begin{lem} The link $L= k \cup \g_1\cup \cdots\cup \g_m$ is prime.
\end{lem} 

\begin{proof} Let $S$ be a $2$-sphere that meets $L$ transversely in exactly two points. The two points must belong to the same component of $L$. Suppose first that this component is a staple $\g_i$. Then $S$ bounds a pair of $3$-balls, one of which contains $k$. The other $3$-ball, which we call $A$, contains an arc of $\g_i$, which must be unknotted as $\g_i$ is trivial. 
It is not possible for $A$ to contain another staple $\g_j$, $j\ne i$, since in that case $S$ would split $k\cup \g_j$, thereby contradicting Lemma 3.1. 
Thus the ball $A$ meets $L$ in an unknotted spanning arc. 

To complete the proof, we need to show that if the two points of $S\cap L$ belong to $k$, then $S$ bounds a ball that intersects $L$ in an unknotted spanning arc.  

Suppose first that $S$ is contained in the interior of some $B_i$. Then $S$ bounds a $3$-ball $A \subset B_i$ meeting $t_i$ in a spanning arc of $A$. Since $(B_i, t_i)$ is a trivial tangle, this spanning arc is unknotted. The lift of $S$ to the $2$-fold cover of $B_i\setminus \g_i$ branched over $t_i$ is a $2$-sphere bounding a $3$-ball that projects to $A$, as $V_i\setminus \tilde \g_i$ is irreducible. Thus $\g_i$ is not contained in $A$, and hence $A$ meets $L$ in an unknotted spanning arc. 

If $S$ is not in the interior of any $3$-ball $B_i$, then we can assume that 
$S\cap (\partial B_1\cup \cdots\cup \partial B_m)$ is a finite collection of pairwise disjoint simple closed curves in which $S$ meets $\partial B_1\cup \cdots\cup \partial B_m$ transversely. Our immediate goal is to show that we can move $S$ without disturbing $L$ setwise so that either $S$ is contained in some $B_i$ or else $S\cap (\partial B_1\cup \cdots\cup \partial B_m)=\emptyset$.

Let $\a$ be a component of $S\cap (\partial B_1\cup \cdots\cup \partial B_m)$ that is innermost in $S$. We can assume that $\a\subset \partial B_i$ and that $\a$ bounds a disk $D\subset S$ such that $D\cap B_j=\emptyset$, for $j\ne i$, and either $D\cap k =\emptyset$ or else $D\cap k$ is one of the two points of $S\cap k$. If $D\cap k = \emptyset$, then either $D\subset B_i$ or $D\subset {\rm cl}(S^3\setminus B_i)$. In the first case, $D$ can be moved off $B_i$, as $L$ is not splittable and $t_{i_1}$ and $t_{i_2}$ are not separated by $D$ in 
$B_i\setminus \g_i$. In the second case, $\a$ also bounds a disk $D'\subset \partial B_i$ such that the cardinality $|D'\cap k|$ is $0,1$ or $2$. If $|D'\cap k|=0$, then the sphere $D\cup D'$ bounds a $3$-ball $A$ such that $A\cap L =\emptyset$, since $L$ is unsplittable or equivalently
${\mathbb S}^3\setminus L$ is irreducible, and we can therefore push $D$ into $B_i$ and thereby remove $\a$ without moving $L$. The case $|D'\cap k|=1$ cannot occur, since $D\cap k=\emptyset$. If $|D'\cap k|=2$, then $D\cup D'$ bounds a $3$-ball outside ${\rm Int}(B_i)$ containing an arc of $k$ and perhaps some of the balls $B_j$. This implies, however, that the crossing of $k$ in $B_i$ is nugatory, contradicting minimality of the projection of $k$. Hence $|D'\cap k|=2$ also cannot occur. 

Assume now that $D\cap k$ is one point, and recall that $\partial D =\a \subset \partial B_i$. Then $\a$ bounds a disk $D'\subset \partial B_i$ meeting $k$ in one point. 

If $D\subset B_i$, then $D\cup D'$ bounds a $3$-ball $A\subset B_i$ meeting $k$ in a spanning arc. Since $(B_i, t_i)$ is a trivial tangle, the arc is unknotted. The irreducibility of the $2$-fold cover of $B_i\setminus \g_i$ branched over $t_i$ implies that $\g_i$ is not a subset of $A$. Hence we can isotop $D$ through $A$ to remove $\a$ while keeping $L$ setwise fixed. 

If $D\subset {\rm cl}({\mathbb S}^3\setminus B_i)$, then the fact that $D$ is an innermost disk (with $\partial D =\a$) in $S$ implies that $D\cap B_j=\emptyset$, for all $j\ne i$, and hence $D\cap k$ is a point in the projection plane. Let $A$ dnote the $3$-ball in ${\mathbb S}^3$ with $\partial A=D\cup D'$ and ${\rm Int}(B_i)$ not a subset of $A$. If $A$ contains any $B_j$, $j\ne i$, then we can move $D'$ slightly off $B_i$ while keeping $k$ setwise fixed to obtain a $2$-sphere $D\cup D'$ such that 
$(D\cup D')\cap (\partial B_1 \cup \cdots \cup \partial B_m)=\emptyset$ and such that $D\cup D'$ bounds two $3$-balls each of which contains at least one of the balls $B_1, \ldots, B_m$. As we will see shortly, this cannot occur, and so $D\cap k$ is a point in one of the four planar arcs of $k$ protruding from $B_i$. These arcs are unknotted  by construction, and no staple $\g_j$ or ball $B_j$ is now in $A$. Hence we can push $D$ back into $B_i$ and remove $\a$, again while keeping $L$ setwise fixed.

We can, therefore, assume that either $S$ is contained in some $B_i$ or 
$S\cap (B_1\cup \cdots \cup B_m)=\emptyset$. As we have seen, if $S$ is in some $B_i$, then $S$ bounds a $3$-ball in $B_i$ meeting $L$ in an unknotted spanning arc. So assume that  $S\cap (B_1\cup \cdots \cup B_m) =\emptyset$. Let $A_1$ and $A_2$ be the two $3$-balls bounded by $S$. Since $k$ is prime, one of $A_1$ and $A_2$, say $A_2$, meets $k$ in an unknotted spanning arc $b$ of $A_2$.

Assume that $S$ is in general position with respect to the projection plane $P$ of $L$. Since the general position isotopy of $S$ can be chosen to fix the two points $x_1$ and $x_2$ of $S\cap k$,  we can assume that $S$ meets $P$ in a simple closed curve containing $x_1, x_2$ together with a collection of simple closed curves bounding disks in $S$. Since we can also assume that $S$ meets a tubular neighborhood $N$ of $k$ (see proof of Lemma 3.3) in two disks, the disks in $S$ bounded by the latter curves belong to the handlebody ${\rm cl}({\mathbb S}^3\setminus{\rm cl}[(\cup_{i=1}^m B_i)\cup N])$, and thus the curves themselves can be removed by cut and paste arguments. Hence there is a simple arc $\b\subset P\cap S$ with $\partial \b = \{x_1, x_2\}$ and a subarc $\a$ of $k$ such that $k=(\a \cup \b)\sharp (\b \cup b)$, where $\b\cup b$ is an unknot, and $k$ is ambient isotopic to $\a\cup \b$. Since the projection of $k$ in $P$ has a minimal number of crossings $m$ (equal to the crossing number of $k$), so does $\a\cup \b$, and so $A_1\supset B_1\cup \cdots \cup B_m \supset \g_1\cup \cdots \cup \g_m$. Therefore, $b\subset P$ and $A_2\cap L = b$. \end{proof}

\begin{lem} The link $L= k \cup \g_1\cup \cdots\cup \g_m$ is hyperbolic. 
\end{lem}

\begin{proof} Let $N$ be a tubular neighborhood of $k$ in $S^3\setminus (\g_1\cup \cdots\cup \g_m)$, and let $N_i$ be a tubular neighborhood of $\g_i$, $i=1, \ldots, m$, such that $N, N_1, \ldots, N_m$ are pairwise disjoint and $N_i \subset {\rm Int}(B_i)$, for each $i$. We also assume  that $N\cap \partial B_i$ is a collection of four meridianal disks of $N$, for each $i$. Set
${\rm Ext}(L)= {\rm cl}({\mathbb S^3}\setminus (N\cup N_1 \cup \cdots \cup N_m))$. With $\tilde \g_i=\tilde  \g_{i_1}\cup \tilde \g_{i_2}\ (i=1, \ldots, m)$, the trivial link
$\g_1\cup \cdots \cup \g_m$ lifts to a $2m$-component link in the $2$-fold cover $M_2$ of $k$, and each $N_i$ lifts to a pair of tubular neighborhoods,
$\tilde N_{i_1}$ and $\tilde N_{i_2}$, of $\tilde \g_{i_1}$ and $\tilde \g_{i_2}$, respectively, in $M_2$. Clearly, $\tilde N_{i_1}\cap\tilde N_{i_2}=\emptyset$
and $\tilde N_{i_1}\cup\tilde N_{i_2}$ is contained in the $2$-fold cover of $B_i$ branched over $t_i$, which is in $M_2$. We set $M= {\rm Ext}(\tilde\g_1\cup \cdots\cup \tilde \g_m)=$ ${\rm cl}(M_2 \setminus \cup_{i=1}^m (\tilde N_{i_1}\cup \tilde N_{i_2}))$, which can be shown to be irreducible by a straightforward application of Lemma 3.2 and the ${\mathbb Z}_2$ sphere theorem \cite{KT}.  Since each of ${\rm Ext}(L)$ and $M$ is an irreducible (in fact, a Haken) $3$-manifold that has torus boundary components and is not a solid torus, it is a standard fact that each of them has incompressible boundary. 

To see that $L$ is hyperbolic, we need to show that $S^3\setminus L$ is not a Seifert fibered space and that every incompressible torus is ${\rm Ext}(L)$ is boundary parallel \cite{thurston2}. That ${\mathbb S}^3\setminus L$ is not Seifert fibered follows
from \cite{BZ}, which yields a geometric description of the unsplittable links in  ${\mathbb S}^3$ with Seifert fibered complements. Each component of such a link can be chosen to be a fiber of some Seifert fibration of ${\mathbb S}^3$. In particular, our link $L$ has four or more components, so if ${\mathbb S}^3\setminus L$ is Seifert fibered, then either (1) each component of $L$ is unknotted; or (2) one or two components are unknotted and each of the remaining components is a nontrivial torus knot (of a given fixed type $(\a, \b)$); or (3) all components are nontrivial torus knots of the same type. Since $L$ has exactly one knotted component but three or more unknotted components, it follows that ${\mathbb S}^3$ is not Seifert fibered. 

We show now that ${\rm Ext}(L)$ is atoroidal, by which we mean that every incompressible torus in ${\rm Ext}(L)$ is boundary parallel. (Our  argument was suggested by that of Case 3 in the proof of Theorem 2 of \cite{GL}.) Suppose first that a torus $T\subset {\rm Ext}(L)$ is incompressible but not boundary parallel and that $T\subset {\rm Int}[B_i \setminus(t_i\cup \g_i)]$, for some $i$. Then the lift $\tilde T$ of $T$ to $V_i \setminus \tilde \g_i$ is either one or two tori. Since $V_i\setminus \tilde \g_i$ is hyperbolic (and thus atoroidal), there is a compressing disk $\tilde D$ for $\tilde T$ in $V_i\setminus \tilde \g_i$ such that 
$g(\tilde D)\cap \tilde D= \emptyset$, or $g(\tilde D) = \tilde D$ and $\tilde D$ meets the fixed point set $\tilde t_i$ of the involution $g$ transversely in a single point \cite{KT} (see also Theorem 3 of \cite{GL}). Let $D$ denote the image of $\tilde D$ under the projection map $V_i\setminus \tilde \g_i \to B_i\setminus \g_i$. If $g(\tilde D)\cap \tilde D = \emptyset$, then the disk $D$ compresses $T$ in ${\rm Int}[B_i\setminus (t_i\cup \g_i)]$, which is a contradiction. If, however, $g(\tilde D)=\tilde D$, then the disk $D$ meets $t_{i_1}$ or $t_{i_2}$ -- say $t_{i_1}$-- transversely in  a single point. We then split $T$ along $D$ to obtain a $2$-sphere $S$ meeting $t_{i_1}$ in two points. As was shown in the  proof of Lemma 3.2, $S$ bounds a $3$-ball $A$ in $B_i\setminus \g_i$ meeting $t_{i_1}$ in  a spanning arc of $A$. It is now clear that $T$ itself must bound the exterior of a nontrivial knot in $B_i\setminus \g_i$, since $T$ is incompressible. This, however, implies that $t_{i_1}$ is a knotted arc, which is a contradiction. Hence $T$ is not contained in ${\rm Int}[B_i \setminus (t_i\cup \g_i)]$, for any $i$. 

On the other hand, the incompressible torus $T\subset {\rm Ext}(L)$ is also not in ${\rm cl}[{\mathbb S}^3 \setminus(N\cup \cup_{i=1}^m B_i]$, as this is clearly a handlebody $(\ne {\mathbb S}^1 \times D^2)$. 

Thus we can assume that $T\cap (\partial B_1\cup \cdots \cup \partial B_m)$ is a finite collection of disjoint simple closed curves along which $T$ and $\partial B_1 \cup \cdots \cup \partial B_m$ meet transversely. Let $\a$ be one of these curves, on $B_i$ say. 

If $\a$ is homotopically trivial on $T$, then it bounds a disk $D\subset T$, and we can assume that $\a$ is innermost on $T$ in the sense that there is no curve $\a'$ in $T\cap (\partial B_1 \cup \cdots \cup \partial B_m)$ such that $\a' \subset {\rm Int}(D)$. Note that $\a \cap \partial (t_{i_1}\cup t_{i_2})=\emptyset$ and that $D$ is properly imbedded in $B_i \setminus (t_{i_1}\cup t_{i_2}\cup \g_i)$ or in ${\rm cl}({\mathbb S}^3\setminus B_i)$. 

\medskip 
{\bf Case 1}\qua $D\subset B_i \setminus (t_{i_1}\cup t_{i_2}\cup
\g_i)$\quad In this case, the disk $D$ lifts to a pair of disks
$\tilde D_1$ and $\tilde D_2$ in $V_i\setminus \tilde \g_i$, each of
which is properly imbedded with $\partial \tilde D_j \subset \partial
V_i$ and $\tilde D_j \cap (\tilde \g_i\cup \tilde t_i) =\emptyset$
($j=1,2$ and $i$ fixed). Since moreover $\partial V_i$ is
incompressible in $V_i\setminus \tilde \g_i$, it follows that
$\partial \tilde D_1$ and $\partial \tilde D_2$ (the lifts of $\a$)
bound disks $\tilde D_1'$ and $\tilde D_2'$, respectively, in
$\partial V_i$ such that $\tilde D_j' \cap \partial \tilde
t_i=\emptyset$ ($j=1,2)$. The projection of $\tilde D_1'\cup \tilde
D_2'$ is a disk $D'\subset \partial B_i$ such that $D'\cap \partial
(t_{i_1}\cup t_{i_2})=\emptyset$ and $\partial D'=\a$, and so $D\cup
D'$ bounds a $3$-ball $A$ in $B_i$ such that $A\cap L=\emptyset$,
since ${\mathbb S}^3\setminus L$ is irreducible. Thus we can isotop
$T$ to remove $\a$.

\medskip{\bf Case 2}\qua $D\subset {\rm cl}({\mathbb S}^3\setminus B_i)$\quad The curve $\a$ bounds two disks $D_1, D_2\subset \partial B_i$ such that $D_1 \cap D_2 =\a$ and $D_1 \cup D_2 = \partial B_i$. If each of ${\rm Int}(D_1)$ and ${\rm Int}(D_2)$ contains a point of  $\partial (t_{i_1}\cup t_{i_2})$, then the minimal number of points in either disk is one or two. Since $D$ contains no points of $k$, however, this minimal number clearly must be two, and since
$|\partial(t_{i_1}\cup t_{i_2})|=4$, each of ${\rm Int}(D_1)$ and ${\rm Int}(D_2)$ must therefore contain two points of $k$. Using $D_1$, say, it follows that $D\cup D_1$ is a $2$-sphere meeting $L$ in  two points of $k$. Since $L$ is prime, $D\cup D_1$ bounds a $3$-ball meeting $L$ in an unknotted arc $b$, a subarc of $k$. Considering $B_i$, this implies that either $k$ consists of two components or the crossing of $t_{i_1}$ and $t_{i_2}$ in $B_i$ is nugatory. Since neither of these is possible, one of $D_1$ and $D_2$ must  miss $\partial(t_{i_1}\cup t_{i_2})$, say $D_1$. By irreducibility of ${\mathbb S}^3\setminus L$, it follows that $D\cup D_1$ bounds a $3$-ball $A\subset {\rm cl}({\mathbb S}^3\setminus B_i)$ such that $A\cap L=\emptyset$, and we can move $T$ to eliminate $\a$. 

Application of Cases 1 and 2 can be used to remove all other curves in $T\cap (\partial B_1 \cup \cdots \cup \partial B_m)$ that are homotopically trivial in $T$ without disturbing the remaining curves. We therefore assume now that $T\cap (\partial B_1 \cup \cdots \cup \partial B_m)$ is a collection of homotopically nontrivial curves in $T$, which must of course be parallel. If this collection is empty, then $T$ is either in some $B_i$ or else $T$ is in the handlebody ${\rm cl}[{\mathbb S}^3\setminus (N\cup \cup_{i=1}^m B_i)]$. Clearly then, a pair of curves, $\a_1$ and  $\a_2$, in $T\cap (\partial B_1 \cup \cdots \cup \partial B_m)$ must bound an annulus $F$ in $T$ with $F\subset B_i\setminus (t_i \cup \g_i)$, for some $i$, and no $\a'$ in $T\cap (\partial B_1 \cup \cdots \cup \partial B_m)$ is contained in ${\rm Int}(F)$. We now show that either $F$ bounds a tubular neighborhood  of $t_{i_1}$or $t_{i_2}$ in $B_i\setminus \g_i$ or else $F$ can be slightly isotoped off $B_i$. 

The curves $\a_1$ and $\a_2$ bound disjoint disks $D_1$ and $D_2$, respectively, in $\partial B_i$, and ${\rm Int}(D_j)\cap \partial(t_{i_1}\cup t_{i_2})\ne \emptyset$ $(j=1,2)$. Since $|{\rm lk}(k,\a_1)|=|{\rm lk}(k,\a_2)|$ and 
$|\partial B_i\cap \partial(t_{i_1}\cup t_{i_2})|=4$, there are three possible cases, two of which we combine into Case (b).

\medskip
{\bf Case (a)}\qua $|{\rm Int}(D_j)\cap \partial(t_{i_1}\cup t_{i_2})|=1$ ($j=1,2)$\quad Since $D_1\cup F\cup D_2$ is a $2$-sphere $S$, it is clear that each of ${\rm Int}(D_1)$ and ${\rm Int}(D_2)$ contains an endpoint of the same arc $t_{i_1}$, say. Isotoping $S$ into $B_i$, it follows that $S$ bounds a $3$-ball $A$ in $B_i\setminus \g_i$ meeting $t_{i_1}$ in an unknotted spanning arc of $A$ (as in the proof of Lemma 3.2.) Isotoping $S$ back to its original position, it follows that $F$ is boundary parallel. (Recall that we began with the original assumption that $T \subset {\rm Ext}(L).)$ 

\medskip
{\bf Case (b)}\qua Either $|{\rm Int}(D_j)\cap \partial(t_{i_1}\cup t_{i_2})|=2$ ($j=1,2)$, or $|{\rm Int}(D_1)\cap \partial(t_{i_1}\cup t_{i_2})|=1$  and 
$|{\rm Int}(D_2)\cap \partial(t_{i_1}\cup t_{i_2})|=3$. (In the second possiblity, the disks' numbering can be switched.) 

Let $F'$ denote the annulus ${\rm cl}[\partial B_i\setminus (D_1\cup D_2)]$, and isotop the torus $F\cup F'$ slightly into ${\rm Int}[B_i\setminus(t_i \cup \g_i)]$ without moving $L$ setwise. As we have seen, the image torus must be compressible in ${\rm Int}[B_i\setminus(t_i \cup \g_i)]$. Now there exist knot exteriors $A_1$ and $A_2$ (at least one of which is a solid torus) such that  ${\mathbb S}^3=A_1\cup A_2$ with $A_1\cap A_2 =F\cup F'$.  One of $A_1$ and $A_2$ (say $A_1$) is in ${\rm Int}[B_i\setminus(t_i \cup \g_i)]$; suppose that $A_1$ is the exterior of a nontrivial knot $k'$, that is, suppose that $A_1$ is not a solid torus. Then the compressing disk $D$ of $F\cup F'$ in ${\rm Int}[B_i\setminus(t_i \cup \g_i)]$ is properly imbedded in $A_2$. The boundary $\partial D$ is not parallel to $\a_1$ (or to $\a_2$) in $F\cup F'$, since each of $\a_1$ and $\a_2$ represents a nontrivial element of $\pi L$ (see Case 2). If $(\partial D, \ell')$ is a meridian-longitude pair for $k'$ (with
$\{\partial D, \ell'\}\subset \partial A_1 = F\cup F')$, it follows that $\a_1$ represents an element of $\pi k' (=\pi A_1)$ of the form $(\partial D)^p(\ell')^q$, where $p,q \in {\mathbb Z}$ with $q\ne 0$. This means, however, that as a simple closed curve in ${\mathbb S}^3$, $\a_1$ must be knotted. But $\a_1$ bounds a compressing disk for $T$ in ${\mathbb S}^3$, and we have a contradiction. Hence $A_1$ is a solid torus. Moving $F\cup F'$ back to its original position, we can thus isotop $F$ through $A_1$ off $B_i$ without disturbing $L$, since $\a_1$ and $\a_2$ are unknotted in ${\mathbb S}^3$. 

Applying Cases (a) and (b) to $T\cup (B_1\cup \cdots\cup B_m)$, we can assume that $T\cup (B_1\cup \cdots\cup B_m)$ is empty except when Case (a) holds for some collection $B_{i_1}, \ldots, B_{i_r}$ ($1\le r\le m$). If
$T\cup (B_1\cup \cdots\cup B_m)= \emptyset$, then $T$ is in the handlebody
${\rm cl}[{\mathbb S}^3\setminus(N\cup \cup_{i=1}^m B_i)]$, which is a contradiction, since $T$ is incompressible in ${\rm Ext}(L)$. Thus we assume that, for some $i$, $T$ meets $B_i$ in an annulus $F$ that is boundary parallel (in $B_i)$ to $\partial N$. The following proposition will enable us to conclude the proof of the lemma.  

\begin{prop} Let $\b = \b_1  \cup \cdots  \cup \b_n$ be a prime 
link in ${\mathbb S}^3$ of $n$ components, and let $T$ be a torus imbedded in ${\mathbb S}^3\setminus \b$. Suppose that $D$ is a compressing disk for $T$ (in ${\mathbb S}^3$) meeting $\b$ transversely in a single point. Then either $\b$  is contained in one component of ${\mathbb S}^3\setminus T$ or else $T$ bounds a tubular neighborhood of $\b_i$, for some $i$. 
\end{prop}

\begin{proof} Assume that $D\cap \b = D\cap \b_1$ is the single point of 
transverse intersection. Assume also that $\b$ is not contained in one
component of ${\mathbb S}^3\setminus T$. If some of $\b_2,\ldots,
\b_n$ are contained in each component, then we surger $T$ along $D$ to
obtain a splitting $2$-sphere $S$ for $\b$ (\figref{hck4}), contradicting
primality.

\begin{figure}[ht!]\anchor{hck4}
\begin{center}
\includegraphics[width=3.3 in]{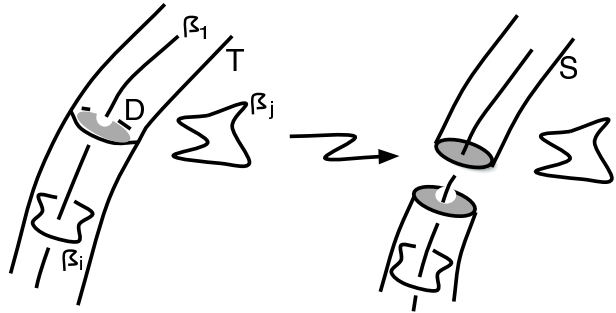}
\caption{Splitting $2$-sphere $S$}
\label{hck4}
\end{center}
\end{figure}

Assume now that $\b_2\cup \cdots \cup \b_n$ lies in the component of ${\mathbb S}^3\setminus T$ not containing $\b_1$ (\figref{hck5}(a)). As in the previous case, surger $T$ along $D$ to obtain a $2$-sphere $S$ (\figref{hck5}(b)). Let $B$ be the $3$-ball with boundary $S$ that does not contain $\b_2\cup\cdots
\cup \b_n$. By primality of $\b$, the $1$-tangle $(B, B\cup \b_1)$ must be trivial. Regard the neighborhood of $D$ removed in surgery as a $1$-handle $h$ with core equal to the part of $\b_1$ not contained in $B$. It is easy to arrange for $h$ to miss $\b_2\cup \cdots \cup \b_n$, since the disk $D$ does not intersect it. Now $B\cup h$ is a solid torus $V$ bounded by $T$. Moreover, the product structure on $h$ extends over $B$ so that $\b_1$ is the core of $V$ (\figref{hck6}). \end{proof}

\begin{figure}[ht!]\anchor{hck5}
\begin{center}
\includegraphics[width=3.3 in]{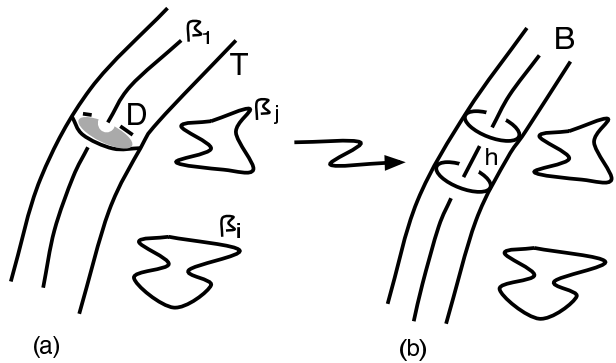}
\caption{Surgery on $T$}
\label{hck5}
\end{center}
\end{figure}

\begin{figure}[ht!]\anchor{hck6}
\begin{center}
\includegraphics[width=3.6 in]{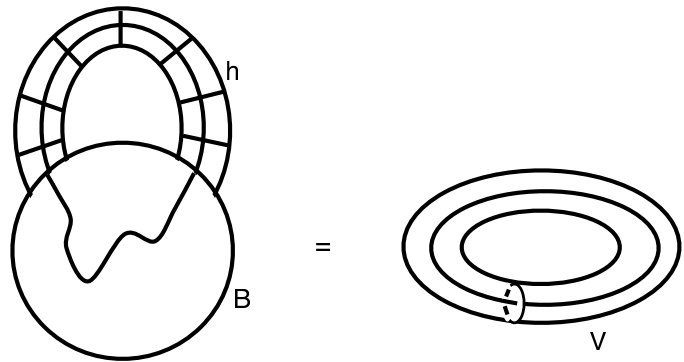}
\caption{$B\cup h$ seen as solid torus}
\label{hck6}
\end{center}
\end{figure}

Continuing with the proof of Lemma 3.3, we have $T\cap B_i=F$, which is boundary parallel to the tubular neighborhood $N$ of $k$. The boundary $\partial F$ is a pair of unknotted curves, $\a_1$ and $\a_2$, bounding disks
$D_1$ and $D_2$ in $\partial B_i$, which are compressing disks for $T$, each meeting $k$ transversely in one point. If $T\cap B_j=\emptyset$, for some $j\ne i$, then $B_j$ is contained in a component $U_1$ of ${\mathbb S}^3\setminus T$. Hence $k\cup \g_j\subset U_1$, and by Proposition 3.4, 
$L\subset U_1$. But if $U_2$  denotes the other component of ${\mathbb S}^3\setminus T$, it is clear that $B_i \cap U_2\ne \emptyset$ and, moreover, that $\g_i \subset U_2$. Thus $T\cap B_j \ne \emptyset$, for all $j$, and $T$ is boundary parallel. Therefore ${\rm Ext}(L)$ is atoroidal, and the proof of Lemma 3.3 is complete. \end{proof}

Since $\g_1$ is unknotted in ${\mathbb S}^3$ and represents the trivial element in $\pi k$, a $1/q_1$-surgery on $\g_1$ changes $k$ into a knot $k_1$ such that $k_1\succeq_1 k$. Now, $\g_2 \subset B_2$, and the $1/q_1$-surgery on $\g_1$ can be regarded as a $(-q_1)$-twist on a disk $D_1\subset B_1$ that is transverse to $k$ such that $\partial D_1 = \g_1$ and $D_1\cap k$ is a set of four points. Thus since $B_1\cap B_2=\emptyset$, it follows that $\g_2$ represents the trivial element of $\pi k_1$, and hence that a $1/q_2$-surgery on $\g_2$ changes $k_1$ into a knot $k_2$ such that $k_2\succeq_1 k_1$. Continuing this process, we arrive at the $m$th stage, in which we do $1/q_m$-surgery on $\g_m$. This changes $k_{m-1}$ into a knot $k_m$ such that $k_m \succeq_1 k_{m-1}$. Thus
$$ k_m \succeq_1 k_{m-1} \succeq_1 \cdots \succeq_1 k_1 \succeq_1 k,$$
and so $k_m \succeq_1 k$. By Thurston's hyperbolic surgery theorem \cite{thurston1}, excluding all but a finite number of possible values of $q_i \in {\mathbb Z}$ for each $i$ assures that $k_m$ is hyperbolic. Hence statement (i) of Theorem 2.2 is proved.

\begin{figure}[ht!]\anchor{ribbondisk}
\begin{center}
\includegraphics[width=2.2 in]{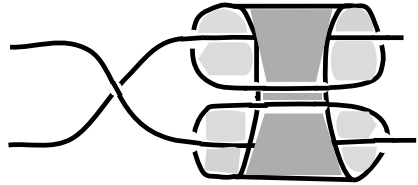}
\caption{Ribbon disk bounded by staple}
\label{ribbondisk}
\end{center}
\end{figure}

In order to prove statement (ii) we observe that the staples $\g_i$ bound pairwise disjoint ribbon disks in the complement of $k$  (\figref{ribbondisk}). The disks can be lifted to the infinite cyclic cover of $k$, and since any two lifts meet only in ribbon singularities, it follows that each $\g_i$ represents an element of the second commutator subgroup of $\pi k$. Hence $1/q$-surgery on $\g_i$ will not change the Alexander invariant
(see Lemma 2 of \cite{nakanishi1}).

Next we prove statement (iii). Let $k_0$ be a hyperbolic knot with trivial Alexander polynomial. Consider the connected sum $k' = k\sharp k_0\sharp \cdots \sharp k_0$ of $k$ with $N$ copies of $k_0$, where $N$ is an arbitrary positive number.  By \cite{SW} there exists a prime knot $k''$ such that $k''\succeq_1 k'$. A proper degree-1 map can be constructed from ${\rm Ext}(k'')$ to ${\rm Ext}(k')$, and hence by \cite{gromov} the simplicial volume of $k''$ is no less than the simplicial volume of $k'$. However, the simplicial volume of $k'$ is at least $N$ times that of $k_0$, which is greater than zero.
Consequently, the simplicial volume of $k''$ can be made arbitrarily large by choosing $N$ sufficiently large.  By part (i) of Theorem 2.2, we can find a hyperbolic knot $\tilde k$ such that 
$\tilde k \succeq_1 k''$. As before, the simplicial volume of $\tilde k$ is at least as large as that of $k''$, and hence the hyperbolic volume of $\tilde k$ can be made arbitrarily large.  

By \cite{SW} and part (ii) of Theorem 2.2, the knots $k', k''$ and $\tilde k$ have the same Alexander invariants. Since $k$ and $k'$ have isomorphic Alexander invariants, so do $\tilde k$ and $k$. 

\begin{figure}[ht!]\anchor{ribbon1}
\begin{center}
\includegraphics[width=4.2 in]{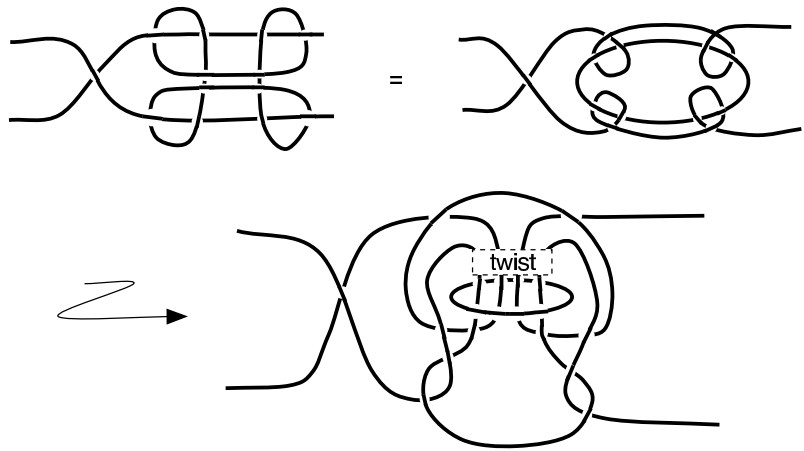}
\caption{Twisting about the staple}
\label{ribbon1}
\end{center}
\end{figure}

\begin{figure}[ht!]\anchor{ribbon2}
\begin{center}
\includegraphics[width=5 in]{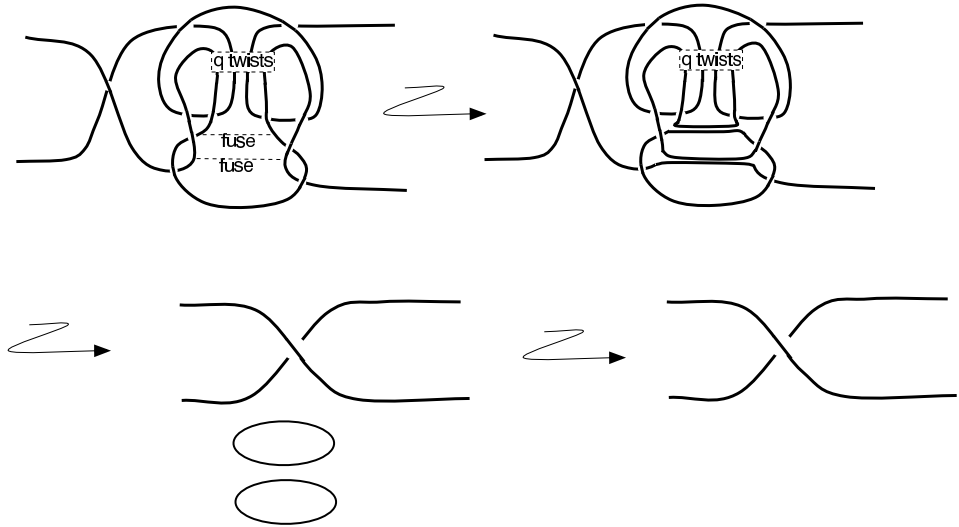}
\caption{Ribbon fusions recovering $k$}
\label{ribbon2}
\end{center}
\end{figure}

Finally we prove statement (iv). The key idea is that $1/q$-surgery on
any staple $\gamma$ converts any knot $k$ to a knot that is ribbon
concordant to $k$. This is immediately seen in Figures \fref{ribbon1}
and \fref{ribbon2}. In \figref{ribbon1}, we see the staple redrawn so
that it bounds an obvious 2-disk. We perform $1/q$-surgery by cutting,
twisting $-q$ full times and reconnecting the strands of $k$ that pass
through the disk.

\figref{ribbon2} shows how a pair of fusions produces two unknotted, unlinked circles that can be shrunk to points. Hence the knot produced from $k$ by surgery is ribbon concordant to $k$.

Recall that we began the proof of Theorem 2.2 by appealing to the main result of \cite{SW}. There we began with any knot $k$, and produced a prime knot by surgery on an unknot $C$ that is not a staple. We complete the proof of Theorem 2.2 (iv) by showing that in fact $C$ can be taken to be a staple. 

\begin{figure}[ht!]\anchor{link}
\begin{center}
\includegraphics[width=1.5 in]{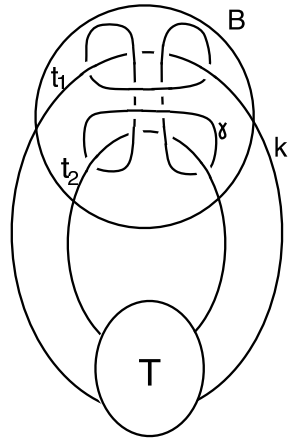}
\caption{$2$-component link $L=k \cup \gamma$}
\label{link}
\end{center}
\end{figure}

According to Proposition 2.5 of \cite{EKT}, we can consider $k$  as the numerator closure $T^N$ of a tangle $T$ that is either prime or rational. Form the 2-component link $L=k \cup \gamma$ (\figref{link}).

Let $(B, t, \gamma)$ be any tangle, where $B$ is a $3$-ball, $t$ is a finite collection of disjoint, properly embedded spanning arcs of $B$, and $\gamma$ is a finite collection of disjoint simple closed curves in ${\rm Int}(B\setminus t)$ such that $t \ne \emptyset$. Following \cite{nakanishi3} and \cite{bleiler}, we will say that $(B, t, \gamma)$ is {\it prime} if it has the following properties.

(i)\qua (No connected summand)\qua Each $2$-sphere in $B$ intersecting $t\cup \gamma$ transversely in two points bounds a $3$-ball in $B$ that meets $t\cup \gamma$ in an unknotted spanning arc. 

(ii)\qua (Disk inseparable\qua) No properly embedded disk in $B\setminus (t\cup \gamma)$ separates $t\cup \gamma$. 

(iii)\qua (Indivisible)\qua Any properly embedded disk $D$ in $B$ such
that $D\cap \gamma =\emptyset$ and such that $D$ meets exactly one
component of $t$ transversely in a single point divides $(B, t,
\gamma)$ into two tangles $(B_1, t', \emptyset)$ and $(B_2, t'',
\gamma)$ such that $t'$ has only one component and that component is
unknotted.

\begin{lem} The tangle $(B, t, \gamma)$ in \figref{link} is prime, where $t=t_1 \cup t_2$. 
\end{lem}

\begin{proof} Form the denominator closure $B^D$. According to the program Snap, a computer program developed at Melbourne University for studying arithmetic invariants of hyperbolic $3$-manifolds (http://www.ms.unimelb.edu.au/~snap/), $B^D$ is a hyperbolic link. Hence $(B, t, \gamma)$ has no connected summand since otherwise $B^D$ would have a connected summand. 

Furthermore, $(B, t, \gamma)$ is disk inseparable since the $2$-fold cover $V\setminus \tilde \gamma$ of $B\setminus \gamma$ branched over $t$ is hyperbolic. A properly embedded disk in $B\setminus (t\cup \gamma)$  lifts to two disks in $V \setminus (\tilde t \cup \tilde \gamma)$, each of which forms a $2$-sphere with a corresponding disk in $\partial V$ that bounds a $3$-ball in $V$ missing $\tilde t \cup \tilde \gamma$. Each of these balls projects to the same $3$-ball in $B \setminus (t\cup \gamma)$.

According to Proposition 1.5 of \cite{nakanishi2}, any tangle that has no connected summand, is disk inseparable, and has at most two spanning arcs is prime. Hence $(B, t, \gamma)$ is prime. \end{proof} 

\begin{lem} The  link $L=k \cup \gamma$ is prime. \end{lem}

\begin{figure}[ht!]\anchor{T_1}
\begin{center}
\includegraphics[width=1.3 in]{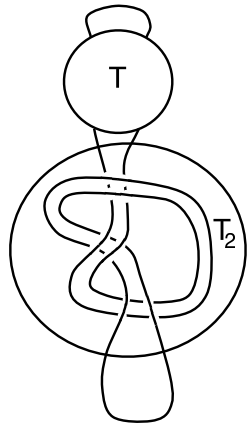}
\caption{The knot $k$ as the numerator closure of $T_1$}
\label{T_1}
\end{center}
\end{figure}

\begin{proof} Since $(B, t, \gamma)$ is prime, this follows immediately from Theorem 1.10 of  \cite{nakanishi2} if $T$ is a prime tangle. 
If $T$ is rational, then we can replace it with a prime tangle $T_1$ such that $T_1^N=k$. The tangle $T_1$ is obtained as a partial sum of $T$ with the prime tangle $T_2$ as shown in \figref{T_1}. It follows from Theorem 3 of 
\cite{lickorish} that $T_1$ is prime, since $T_2$ is prime. Hence again $L$ is prime. \end{proof}

The remaining argument of \cite{SW} applies now, completing the proof of Theorem 2.2 (iv).  

\begin{proof}[Proof of Corollary 2.3] 
Let $k_0$ be the untwisted double of a trefoil. Corollary 5 and
Theorem 1 of \cite{livingston} together imply that the $4$-ball genus
of the connected sum $k\sharp k_0\sharp \cdots \sharp k_0$ can be made
arbitrarily large by increasing the number of summands $k_0$. (The
results of \cite{livingston} are convenient for us, but earlier work
of Rudolph \cite{rudolph} could be used instead.) We replace $k$ by
$k\sharp k_0\sharp \cdots \sharp k_0$, which has the same Alexander
invariant, and apply Theorem 2.2. Since the resulting knot $\tilde k$
is (ribbon) concordant, the two knots have the same $4$-ball
genus. \end{proof}

\bibliographystyle{gtart}

\Addresses\recd
\end{document}